\documentclass[a4paper,12pt]{amsart}
\usepackage{amsfonts, amsmath,amscd,amssymb}
\newtheorem{theorem}{Theorem}[section]

\newtheorem{prop}[theorem]{Proposition}
\newtheorem{lem}[theorem]{Lemma}
\theoremstyle{definition}
\newtheorem{defn}[theorem]{Definition}

\theoremstyle{remark}
\newtheorem{rem}[theorem]{Remark}
\numberwithin{equation}{section}
\def\prf{\textsc{Proof}. }
\def\endprf{\hfill$\Box$\vskip .25cm}

\def\Cl{\mathop{\operatorname{\mathrm Cl}}\nolimits}
\def\Id{\mathop{\operatorname{\mathrm Id}}\nolimits}
\def\ind{\mathop{\operatorname{\mathrm ind}}\nolimits}
\def\TVect{\mathop{\operatorname{\mathrm Vect}}\nolimits}
\def\SO{\mathop{\operatorname{\mathrm SO}}\nolimits}
\def\Spin{\mathop{\operatorname{\mathrm Spin}}\nolimits}
\title{On the Twisted KK-Theory and Positive Scalar Curvature Problem}
\author{Do Ngoc Diep}
\address{
Institute of Mathematics, Vietnam Academy of Science and Technology, 18 Hoang Quoc Viet Road, Cau Giay District, 10307, Hanoi, Vietnam,
and\newline c/o Vietnam Institute for Advanced Study in Mathematics, 1 Ta Quang Buu, Hanoi, Vietnam}
\email{dndiep@math.ac.vn}

\thanks{\bf The author's thanks are addressed to VIASM for a financial support for a scientific stay from December 07, 2015 to January 06, 2016, at that Intsitution.}
\subjclass{Primary: 19L50 ; Secondary: 19K56}
\date{\textbf{\today}}
\begin{document}
\maketitle

\maketitle
\begin{abstract}
Positiveness of scalar curvature and Ricci curvature requires  vanishing the obstruction $\theta(M)$ which is computed in some KK-theory of C*-algebras index as a pairing of spin Dirac operator and Mishchenko bundle associated to the manifold. U. Pennig had proved that the obstruction $\theta(M)$ does not vanish if $M$ is an enlargeable closed oriented  smooth manifold of even dimension larger than or equals to 3, the universal cover of which admits a spin structure.  Using the equivariant cohomology of holonomy groupoids we prove the theorem in the general case without restriction of evenness of dimension. Our groupoid method is different from the method used by B. Hanke and T. Schick in reduction to the case of even dimension.
\vskip .1cm
\textit{Key terms}: scalar curvature; Ricci curvature; KK-theory; index theory; 

\textbf{2010 Mathematics Subject Classification}: \textit{Primary}: 19L50; \textit{Secondary}: 19K56

\end{abstract}
\tableofcontents

\section{Introduction}

The problem of existence of positive scalar curvature on a smooth manifold $M$ is reduced to computation of the obstruction $\theta(M)$, see \cite{rosenberg}. It is evaluated in the KO-theory $KO_n(C^*_\mathbb{R}(\pi_1(M)))$ of group C*-algebras of the fundamental group $\pi_1(M)$. The Gromov, Lawson and Rosenberg conjecture states that this invariant is adecquate to existence of positive scalar curvature metric on $M$, if $\dim M \geq 5$. The conjecture was proven by Stolz \cite{stolz} for the case when $M$, or at least his universal covering $\tilde{M}$ is a spin manifold. He constructed $\theta(M)$ to take values in the $KO$ group $KO(C^*_\mathbb{R}(\hat{\pi}\to \pi_1(M))$ of the group C*-algebra of the twisted fundamental group $\hat{\pi}$ as the central two-fold extension of the fundamental group $\pi_1(M)$. The element $\theta(M)$ can be expressed as the pairing 
$$\theta(M) = \ind(D^\nu_+) =  [\nu^S] \otimes_{C(M,\Cl(M))} [D^S]$$
of the twisted Dirac class $[D^S]$ from the K-homology $KK(C(M, \Cl(M)), \mathbb C)$ and the twisted Mishchenko class $[\nu^S]$ from the K-cohomology $KK(\mathbb C, C(M, \Cl(M)\otimes C^*(\hat{\pi}\to \pi(M)))$.  He restricted to the case where the $\pi_1(M)$ modules $E$ are of finite dimension. We remark that the construction of the cup product of pairing between the Dirac operator for  Clifford module and the Mishchenko bundle 
is also generalized  for the infinite dimensional finitely generated ones. 

Our approach is to consider the two-fold covering 
$$\CD 0 @>>>\mathbb Z/2\mathbb Z @>>>\hat{G} @>>> G @>>> 1 \endCD $$ 
of the holonomy groupoid $G\to M\times M$  which is constructed following the commutative diagram
 
$$\CD 
@. 1 @= 1 @. @. \\
@. @VVV @VVV @. @. \\
0 @>>> \mathbb Z/2\mathbb Z @= \mathbb Z/2\mathbb Z @>>> 0 \\
@. @VVV  @VVV  @VVV \\
1 @>>> \hat{\pi} @>>> \hat{G} @> s,t>> M\times M@>>> pt \\
@. @VVV   @VVV  @| @.\\
1 @>>>\pi_1(M) @>>> G @> s,t>> M\times M @>>>pt\\
@. @VVV @VVV @VVV @. \\
@. 1 @= 1 @>>>0 @. \\
\endCD$$

Let us describe in more details the situation. Let $\tilde{M}$ be the universal cover of $M$. The groupoid of holonomy $G$ is consisting of all homotopy classes of free-end curves on $M$ 
$$G =  \left\{(x,y,[\gamma]) \left| x,y \in M, [\gamma] ={ \mbox{homotopy class of curve }\gamma; \atop s(\gamma)=\gamma(0) = x, r(\gamma)=\gamma(1) = y } \right. \right\}$$ The source- and target- maps $s,t: P \to M$ provide the corresponding groupoid structure on $P$.  The product of two free-end curves $\gamma_1,\gamma_2$ can be defined if $t(\gamma_1) = y_1 = x_2 = s(\gamma_2)$, i.e.
$$(x_1,y_1,[\gamma]) \circ (x_2,y_2,[\gamma_2]) = (x_1,y_2,[\gamma_1\circ\gamma_2]).$$
The C*-algebra $C^*(G)$ is defined in a standard way: convolution of two absolutely integrable function on $P$ is defined as
$$(f *g)(\gamma) = \int_{\gamma_1\circ\gamma_2 = \gamma}f(\gamma_1)g(\gamma_2) d\gamma_1.$$ The C*-algebra $C^*(G)$ is defined as the completion of the Banach algebra $L^1(G)$ with respect to $d\gamma$ over the cover $\tilde{M}$ of $M$. From the results of theory of groupoid. 

\begin{theorem}[Main Theorem]
Every closed smooth manifold the universal covering of which has a spin structure and which is  enlargeable has the nonvanishing theta invariant, $\theta(M) \ne 0$. 
\end{theorem}
\begin{rem}
We show that the assumptions of even dimension in the work is removed, because we pass to the equivariant theory of holonomy groupoid $\CD G @>s,t>> M\times M\endCD$.  
\end{rem}

\begin{rem}
Combine our result with the well-known result of Gromov Lawson: ``every simply connected closed non-spin manifold of dimension $n \geq 5$ admits a positive scalar curvature meric" one deduces that for the enlargeable simply connected closed manifolds of dimension $n\geq 5$ the spin-structure is decisable condition for the existence of positive scalar curvature metric: those manifolds admit a spin-structure if and only iff the obstruction to existence of a PSC metric vanishes.

\end{rem}
The paper is organised as follows: In \S2 we introduce some preparation and in particular introduce the notion of bundle gerbe, Mishchenko bundle gerbe and Mischenko bundle. Section 3 is devoted to proof of the main theorem. In \S4 we related the problem with the Lie algebroid cohomology.

\begin{rem} 
In the paper (\cite{hankeschick}, Proof of Thm. 4.2) B. Hanke and T. Schick reduced the general case to the case that M has even dimension, by using the diagram composing of homotopy classification of K-homology and the Baum-Connes assembley map from the K-homology of classifying spaces to the corresponding K-homology of C*-algebras of the fundamental groups, and their corresponding images underer multiplication with the fundamental calss $[\mathbb S^1]$.  Our method of passing to groupoids, is quite different from those high developed techniques. 
\end{rem}

\section{Preliminaries}
We will use the same notation as in the paper of U. Pennig \cite{pennig}.  Let us remind and modify some notions from that work for our usage. 
 
\subsection{Bundle Gerbes}
\begin{defn}[bundle gerbe]   Consider a surjective submersion $Y \to M$ and a real line bundle $L\to Y^{[2]} = Y \times_M Y$ and its pullback $pr_{ij}^*L$ over $Y^{[2]}$ with respect to the projections $pr_{ij}: Y^{[3]} \to Y^{[2]}$ on factors. It is called a $\mathbb Z/2\mathbb Z$-\textit{bundle gerbe} if the following diagram is commutative
$$\CD (pr_{12}^*L \otimes pr_{23}^*L) \otimes pr_{34}^*L @= . @= pr_{12}^*L \otimes (pr_{23}^*L \otimes pr_{34}^*L) \\
@V{\mu \otimes id }VV @. @VV{id \otimes \mu}V \\
pr_{13}^*L \otimes pr_{34}^*L @>>\mu> pr_{14}^*L @<<\mu< pr_{12}^*L \otimes pr_{24}^* L  \endCD$$ 
The bundle gerbe $\delta Q= pr_1 Q \otimes pr_2 Q \to Y^{[2]}$ which is the pullback from a line bundle $Q \to Y$ under projections $pr_i: Y^{[2]} \to Y$ is called \textit{trivial bundle gerbe}.
\end{defn}

Let us remind from \cite{pennig} that if $M$ is a closed smooth oriented manifold of dimension $\dim(M) \geq 3$, then the oriented frame bundle $P=P_{\SO}$ is a principal bundle,
we have the short exact sequence of spaces
$$\SO(n) \hookrightarrow P_{\SO( M)} \twoheadrightarrow  M $$  
therefore we have a part of long exact sequence
$$\CD \pi_2( M) @>>> \pi_1(\SO(n)) @>>> \pi_1(P_{\SO( M)}) @>>> \pi_1( M) @>>> 1,\endCD$$ then 
 for the universal covering
$$\SO(n) \hookrightarrow P_{\SO(\tilde M)} \twoheadrightarrow \tilde M $$  
we have a similar exact sequence
$$\CD \pi_2(\tilde M) @>>> \pi_1(\SO(n)) @>>> \pi_1(P_{\SO(\tilde M)}) @>>> \pi_1(\tilde M) = 1\endCD$$
We suppose that the universal covering $\tilde{M}$ admids a spin structure. 
If the map $f: \mathbb S^2 \to \tilde M$ is a representative of the second homotopy group $\pi_2(\tilde M)$ then the map $\pi_2(\tilde M) \to \pi_1(\SO(n))$ sends it to a map $\varphi_f : \mathbb S^1 \to \SO(n)$ obtained from the pullback $f^*P_{\SO(\tilde M)}$ and it must be factorized through a map $\mathbf S^1 \to \Spin(n)$ which is nullhomotopic and the same is also the map $\pi_2(M)\to \pi_1(\SO(n))$. One has therefore the short exact sequence
$$\CD 1 @>>> \mathbb Z/2\mathbb Z @>>> \pi_1(P_{\SO}) @>>> \pi_1(M) @>>> 1\endCD$$
and therefore we also have an exact sequence
$$\CD 1 @>>> \mathbb Z/2\mathbb Z @>>> \pi_1(\hat G) @>>> \pi_1(M\times M) @>>> 1.\endCD$$

Following U. Pennig \cite{pennig} we have 
\begin{defn}
The lifting bundle gerve $L_G \to \hat G ^{[2]}$ associated to the central extension $\hat G \to \tilde M \times \tilde M$ is called the \textsl{Mishchenko bundle gerbe}.
\end{defn}

\subsection{Construction of $\theta$-Invariance}

\begin{defn}
The fibered product $\nu_r = \hat G \times_{M\times M} C^*_r(\hat \pi_1(M))$ is called \textsl{twisted Mishchenko bundle}.
\end{defn}

\begin{defn} For any $C^*(\hat\pi)$-module $E$,
the index of the Dirac operator $D^E$ representing a class $[D^E]\in KK(\mathbb C, C(M,\Cl(M)) \otimes C^*(\hat{\pi}))$
$$\theta^E(M) = \ind(D^E_+) = [\nu^E] \otimes_{C(M,\Cl(M))} [D^S]$$ is called the \textsl{$\theta$-invariant} of $M$. In the case of regular representation or $E=C^*(\hat{\pi})$ one denote simply $[\nu^S]$ and $\theta(M)$
\end{defn}

\subsection{Equivariant KK-Theory of Groupoids}
We will show in the next section that
there is a natural isomorphism between equivariant KK-groups of holonomy groupoid and the corresponding KK-group for $M$, and the same for cyclic homologies.

\begin{rem} It reduces the theory of cohomology of the manifold $M$ to the 
equivariant cohomology of holonomy groupoids. Later we will reduce the equivariant cohomology theory of holonomy groupoids to the cohomology of correspoding Lie algebroids, which is computed in more convenient situation of linear algebra.  \end{rem}

\section{Proof of the Main Theorem}
The idea is to apply the main theorem of U. Pennig \cite{pennig}. In order to do so we do check all the condition to apply such a result to the two-fold covering of the holonomy groupoid. 
\begin{lem}
There is a natural two-fold covering $\hat{G}$ of the groupoid of holonomy  $\CD G @>s, t >> M\times M\endCD$.  .
\end{lem}
\prf
The two-fold covering $\hat{G}$ 
$$\CD 0 @>>>\mathbb Z/2\mathbb Z @>>>\hat{G} @>>> G @>>> 1 \endCD $$ 
of the holonomy groupoid $\CD G @>s, t>> M\times M\endCD$  which is constructed following the commutative diagram
 
$$\CD 
@. 1 @= 1 @. @. \\
@. @VVV @VVV @. @. \\
0 @>>> \mathbb Z/2\mathbb Z @= \mathbb Z/2\mathbb Z @>>> 0 \\
@. @VVV  @VVV  @VVV \\
1 @>>> \hat{\pi} @>>> \hat{G} @>r \times s>> M\times M@>>> pt \\
@. @VVV   @VVV  @| @.\\
1 @>>>\pi_1(M) @>>> G @>r \times s>> M\times M @>>>pt\\
@. @VVV @VVV @VVV @. \\
@. 1 @= 1 @>>>0 @. \\
\endCD$$
\endprf

\begin{lem}
If $L$ is a bundle gerbe over $\tilde{M}$ then there is a bundle gerbe $L_G$ over $\hat{G}$. 
\end{lem}
\prf
Indeed, if we have a bundle gerbe $L\to \tilde{M}$ then $L \times L$ is a bundle gerbe over $\tilde{M} \times \tilde{M}$ and therefore we may construct $L_G = (r \times s)^*(\tilde{L} \times \tilde{M})$ as the corresponding pullback. 
\endprf

The following Lemma is trivial.
\begin{lem}
If $\dim M \geq 2$ then $\hat{G}$ is a manifold of even dimension larger  than 3.
\end{lem}
\prf
Locally $\hat{G}$ is homeomorphic to $M \times M$.
\endprf

\begin{lem}
if the manifold $M$ is enlargeable, then $\tilde{G}$ is also enlargeable.
\end{lem}
\prf
If for any positive $\varepsilon$ there exists an $\varepsilon$-contraction map $(M,g)\to (\mathbb S^n,g_0)$, then there exists also an   $\varepsilon$-contraction map $(\hat{G},g)\to (\mathbb S^{2n},g_0)$ because locally $\hat{G}$ and $\hat{M}\times \tilde{M}$ are diffeomorphic to each-another and the same is locally, for $\mathbb S^n \times \mathbb S^n$ and $\mathbb S^{2n}$.
\endprf

\begin{lem}
If the universal covering $\tilde{M}$ of the manifold $M$ has a spin structure then $\hat{G}$ has also a spin structure.
\end{lem}
\prf
If there is a spin structure on $\tilde{M}$ then the same one exists on $\tilde{M}\times \tilde{M}$ and therefore on $\hat{G}$.
\endprf

Let us denote by $P$ the principal bundle of frames over $M$ and $\tilde P$ its two-fold covering.
The next lemma reduces the computation of the $K$-group of $M$ with spin structure on $\tilde{M}$ to the equivariant one on the holonomy groupoid.
\begin{lem}
The equivariant $K$-group $KK^0_{\hat G}(C^*(\hat{G}), C(G,\Cl(G))\cong$ \newline $ KK^0_{\pi_1(M)}(C(\tilde{P}), C(P,\Cl(P))$ is isomorphic to the $K$-group \newline $KK^0(C(\tilde{M}), C(M,Cl(M)).$
\end{lem}
\prf
When we fix a connection on the  holonomy groupoid we have the principal bundle $\hat P= \hat P_{\SO}$. The first isomorphism is natural. The second isomorphism is also obtained from the natural diffeomorphism $\hat P\slash \SO(M) \cong M$ from the exact sequence of spaces
$$\CD 1@>>> \SO(n) @>>> \hat P_{\SO(M)} @>>> M.\endCD$$ 
\endprf

\begin{lem}
The Dirac operator $D^E_G$ on the holonomy groupoid $G$ can be decomposed as the sum of two copies of the corresponding Dirac operators $D^E_M$ on $M \times \{pt\}$ and $\{pt\} \times M$
$$D^E_G = D^E_M \otimes \Id + \Id \otimes D^E_M.$$
\end{lem}
\prf Locally $G$ is homeomorphic to $M \times M$. From the local definition of Dirac operator we have: if the index of $D^E_M$ vanishes, so does the index of $D^E_G$. But the index of the Dirac operator class $[D^E_G]$ does not vanish, following \cite{pennig}. \endprf

\begin{lem}
The equivariant cylic cohomology of the holonomy groupoid is isomorphic with the cyclic cohomology of Lie algebroid. 
\end{lem}
\prf
There is a natural isomorphism of equivariant differential forms and equivariant tangent vector fields. Therefore, 
the natural isomorphism between equivariant cohomology on the holonomy groupoid and the cohomology of Lie algebroids.
\endprf

\section{Lie algebroid cohomology}
The computation of equivariant cyclic cohomology of the groupoid of holonomy is reduced to the same one for Lie algebroids of tangent fields. In many case it should be easear than the original ones.
 
It is the Lie algebroid $\TVect(M)$ of equivariant free-end tangent to curve vector fields. For two free-end tangent vector fields, their Lie bracket can be defined when the sources and rivals are coinc\"ided respectively. Trivial verification show that the Lie brackets, when they are  defined satisfied the axioms of Lie theory of vector field

We do compute it for the cohomology.  
\begin{prop}
The equivariant groupoid cyclic cohomology $HC^*_{\hat{G}}(\hat G; \mathbb R)$ is isomorphic to the cyclic cohomology $HC^*(\mathfrak g; \mathbb R)$ of Lie algebroid $\mathfrak g =\TVect(M)$ of invariant vector fields.
\end{prop}
\prf
When we fix a trivialization of $\hat G$ we come to the oriented frame bundle $\hat P$. For $\hat P$ and corresponding $\mathfrak p$ as Lie subgroupoid and Lie subalgebroid, we have natural isomorphism of equivariant cyclic theories.
\endprf

\section*{Acknowledgments}
The paper was done in the framework of a project of research group ``K-Theory and Noncommutative Geometry'', at VIASM. The author thanks the Institute for a scientific stay support during the 2015 December and 2016 January. 

The sincere thanks are extended to all members of the group, especially Professors A.S. Mishchenko, V. M. Manuilov, G. Sharygin, F. Popelensky and visitors during the seminars for the constructible very useful discussions.

\end{document}